# MODELING AND ANALYSIS METHODS FOR EARLY DETECTION OF LEAKAGE POINTS IN GAS TRANSMISSION SYSTEMS

Ilgar Aliyev

Head of Department, Azerbaijan Architecture and Construction University, Baku, Azerbaijan

## ABSTRACT

*Early detection of leaks in gas transmission systems is crucial for ensuring uninterrupted gas supply, enhancing operational efficiency, and minimizing environmental and economic risks. This study aims to develop an analytical method for accurately identifying leak locations in gas pipelines based on unsteady gas flow dynamics. A novel approach is proposed that utilizes pressure variations at the inlet and outlet points to determine the minimum fixation time ($t = t_1$) required for real-time leak detection. Through mathematical modeling and numerical analysis, the study demonstrates that the ratio of pressure drops at different points along the pipeline can be effectively used to pinpoint leakage locations. The results indicate that the proposed method significantly improves detection accuracy and response time, making it a viable solution for integration into gas pipeline monitoring and control systems.*

## KEYWORDS

*leakage, modeling, optimal, analytical method, operational*

## 1. İNTRODUCTION

The transportation of natural gas through pipelines is a critical component of the global energy infrastructure. These pipeline networks serve as the backbone for supplying energy to industries, businesses, and households. However, ensuring the integrity and safety of gas pipelines remains a major challenge due to the risks associated with gas leaks. Leakage in gas pipelines can lead to severe consequences, including environmental pollution, financial losses, and safety hazards such as fires or explosions. To ensure the efficient and safe operation of complex gas pipelines, it is crucial to promptly detect, eliminate, and minimize natural gas loss by addressing the issues related to the detection and resolution of gas leaks.

Gas pipeline leaks can result from various factors, including corrosion, mechanical damage, material defects, and natural disasters [1]. For underground trunk pipelines, leaks are particularly dangerous because they can go undetected until significant damage has occurred. The presence of microcracks, caused by environmental factors or pipeline aging, can gradually lead to the failure of pipeline sections under pressure. These failures may not only disrupt gas supply but also pose risks to human life and property.

Conventional methods for leak detection, such as visual inspection, odor detection, and ground-based monitoring, often prove inadequate for modern gas pipeline systems. More advanced techniques, including acoustic, hydraulic, and differential pressure methods, have been introduced to improve leak detection efficiency. However, many of these methods require complex sensor networks, high computational power, or extensive calibration, making their implementation costly and sometimes impractical. Consequently, there is a strong need for an optimized, real-time approach that minimizes detection time and operational costs while maintaining high accuracy in leak localization.

This study focuses on developing an analytical model to detect and localize gas leaks in transmission pipelines using unsteady gas flow parameters. Unlike traditional methods [2, 3, 4, 5, 6, 7, 8], this approach leverages pressure variations at the inlet and outlet points of the pipeline to determine leak locations within a minimal time frame. By analyzing transient flow conditions and pressure fluctuations, the proposed method enhances leak detection capabilities without requiring extensive sensor deployment.

Existing techniques for gas leak detection include [9,10,11]:

Pressure wave analysis: Detects sudden pressure changes but requires continuous monitoring and may introduce errors due to wave dissipation.

Hydraulic correlation methods: Estimates leakage points based on pressure measurements but may produce inaccuracies due to flow variations.

Acoustic leak detection: Identifies characteristic ultrasonic noise from gas leaks but is expensive and requires high-density sensor placement.

Differential flow balance methods: Offers high accuracy for detecting small leaks but requires precision flow meters at multiple locations.

Our approach builds upon these methods by introducing a new computational framework that determines leakage locations based on transient flow analysis while optimizing sensor placement to minimize costs. The proposed methodology allows for real-time monitoring, making it suitable for automated pipeline control systems.

Objectives of the Paper: To develop an analytical model for real-time gas leak detection based on transient gas flow analysis;To determine the minimal fixation time (t = $t_1$) required to accurately detect leaks and optimize emergency response measures;To propose a cost-effective solution by minimizing the number of required sensors while maintaining high accuracy in leak localization (Figure 1). The sensors installed at the beginning and the end of each section of the gas pipeline are queried at the control point, which increases the efficiency and performance of the control point.

## 2. Materials and Methods

In this study, a combination of software-based analysis and pressure sensor data was used to develop a model capable of detecting gas leaks with high precision and minimal response time. The approach focuses on optimizing cost efficiency by minimizing the number of required sensors while maintaining accuracy in leak detection.

The fundamental mathematical model for gas flow in pipelines under unsteady conditions is based on the conservation laws of mass, momentum, and energy. A system of partial differential equations was derived and solved numerically to analyze transient gas flow behavior during leakage events. The modeling framework assumes compressible flow and considers the impact of gas pipeline elasticity on pressure wave propagation. To validate the proposed model, numerical simulations were conducted using MATLAB and ANSYS Fluent. The pipeline was segmented into discrete sections, with pressure variations recorded at both inlet and outlet points.

A real-time data processing algorithm was developed to identify leak locations using pressure sensor readings. The algorithm operates as follows:

Data Collection: Pressure readings are continuously gathered from the pipeline.
Leak Localization: Using an optimization-based approach, the algorithm estimates the leak position by analyzing differential pressure variations between inlet and outlet points.
Decision Support: The system provides real-time alerts to operators, enabling prompt action to mitigate gas loss.
By integrating the developed model into an automated monitoring system, gas leaks can be detected with minimal delay, enhancing operational safety and efficiency.
 Advantages of the Proposed Method:
High Sensitivity-Detects even small leaks by analyzing transient pressure variations.
Cost Efficiency: Reduces sensor requirements while maintaining detection accuracy.
Real-Time Monitoring: Enables continuous surveillance and rapid response to leakage incidents.
Scalability: Can be adapted for different pipeline lengths and configurations.
This section provides a comprehensive description of the techniques used for gas leak detection, ensuring a robust and efficient approach for industrial applications.

When preparing solutions for the given problem, it is advisable to study non-stationary processes and obtain mathematical expressions that can manage the occurring physical processes. To achieve this goal, it is necessary to use non-stationary models. In order to strike a balance between adequacy and simplicity, the requirements for mathematical models compel us to closely monitor the physical processes of gas flow through pipes and require the development of specific models for each transition process.

We accept the following technological diagram of the gas pipeline, which reflects the physical processes being observed or studied (Figure 1).

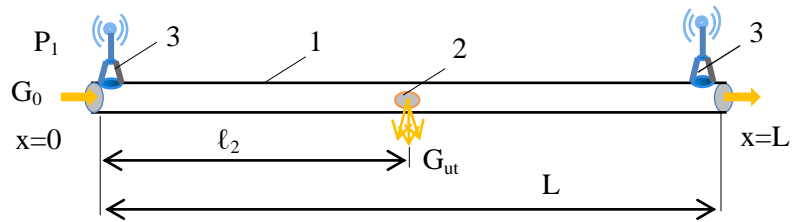

Figure 1. Diagram of a Gas Pipeline with Leakage
1 - Gas pipeline, 2 - Location of the leakage point in the pipeline,
3 - Wireless pressure sensor measuring gas flow pressure, L - Length of the gas pipeline,
$P_1$ - Initial pressure of the gas pipeline, $P_2$ - Final pressure of the gas pipeline,
$G_0$ - Mass flow rate at the beginning of the gas pipeline, $G_{ut}$ - Mass flow rate of the leaking gas,
$\ell_2$ – distance from the starting point of the gas pipeline to the leakage point.

Wireless pressure sensors are intended to monitor pressure fluctuations at the beginning (x=0) and end (x=L) of the gas pipeline. As seen in Figure 1, it is proposed to install these pressure sensors on the pipeline at the starting and ending sections of the gas pipeline. These sensors enable the monitoring of pressure variations in the pipeline through the internet at the control station.

## 3. Analysis of the mathematical model for unsteady gas flow in the pipeline.

For non-stationary gas flow, I.A. Charny presented a system of partial differential equations related to gas flow parameters, such as pressure P(x, t) and gas mass flow G(x, t), as functions of two variables — time t measured along the axis and distance x. Then, for the problem we are considering, the system of partial derivatives of nonlinear equations will be as follows [12,13].

$$\begin{cases} -\dfrac{\partial P}{\partial X} = \lambda \dfrac{\rho V^2}{2d} \\ -\dfrac{1}{c^2}\dfrac{\partial P}{\partial t} = \dfrac{\partial G}{\partial X} + 2aG_{ut}(t)\delta(x-\ell_2) \\ P = \rho z RT \end{cases} \quad (1)$$

Here,
T - average temperature of the gas
z = z(p, T) - compressibility coefficient
c - speed of sound propagation of the gas for an isothermal process, m/second
P - pressure at the ends of the gas pipeline, Pa (Pascal)
x - longitudinal coordinate of the gas pipeline, m (meter)
t - time coordinate, second
λ - hydraulic resistance coefficient
ρ - average density of the gas, kg/ $m^3$ (kilogram per cubic meter)
V - average velocity of the gas flow, m/second
d - diameter of the gas pipeline, m (meter)
G - mass flow rate of the gas
$\delta(x-\ell_2)$ - Dirac function

It is clear that solving the non-linear system of equations (1) is not possible because these types of equations cannot be integrated. Approximate methods are used in engineering calculations. In other words, equations are linearized. The following linearization methods are used to solve the equations: It has been shown during calculations that Chernov's linearization is more favorable when the pipelines are damaged. The expression characterizing this linearization is as follows: [14].

$$2a = \lambda \dfrac{V}{2d} \quad (2)$$

Here, 2a represents the Chernov linearization factor, in which the system of equations (1) and (2) will be as follows:

$$\begin{cases} -\dfrac{\partial P}{\partial X} = 2aG \\ -\dfrac{1}{c^2}\dfrac{\partial P}{\partial t} = \dfrac{\partial G}{\partial X} + 2aG_{ut}(t)\delta(x-\ell_2) \end{cases} \quad (3)$$

Here, G=ρv; P=ρ·$c^2$

If we take the derivative with respect to x from the first equation of the system (3) and determine − from the second equation, then we obtain the equations in the following form.

$$\begin{cases} -\dfrac{\partial^2 P}{\partial x^2} = 2a\dfrac{\partial G}{\partial X} \\ \dfrac{\partial G}{\partial x} = -\dfrac{1}{c^2}\dfrac{\partial P}{\partial t} - 2aG_{ut}(t)\delta(x-\ell_2) \end{cases} \qquad (4)$$

After substituting the expression for $\dfrac{\partial G}{\partial x}$ into the first equation of the system (4), we obtain the heat transfer equation for solving the mathematical problem.

$$\dfrac{\partial^2 P}{\partial X^2} = \dfrac{2a}{c^2}\dfrac{\partial P}{\partial t} + 2aG_{ut}(t)\delta(x-\ell_2) \qquad (5)$$

Therefore, to fully define the given process, solving equations (5) should be sufficient, provided that the initial and boundary conditions are correctly specified. In the initial conditions, the distribution of the sought function (pressure) along the length of the pipe at t=0 is given when inertia forces are not considered. The boundary conditions are provided to ensure the prescribed law-abiding changes in pressure and flow at the endpoints of each segment of the pipe throughout the considered time. Properly specifying the boundary conditions completes the mathematical model of the process and allows for a comprehensive and accurate investigation of the physical events occurring.

Thus, based on the investigated physical process and Figure 1, we accept the expressions representing the initial and boundary conditions as follows.

$$t=0;\ P(x,0)=P_1-2aG_0 x:$$

For the main sections of the gas pipeline considered in the process, we assume the following boundary conditions (where flows are measured at the initial and final points):

At x=0 point $\dfrac{\partial P(x,t)}{\partial x} = -2aG_0(t);$ At x= L point $\dfrac{\partial P(x,t)}{\partial x} = -2aG_s(t)$

Here, $G_0(t)$, $G_s(t)$ və $G_{ut}(t)$ respectively represent the mass flows measured at the beginning, end, and leakage point of the pipeline at different times, $\dfrac{Pa \cdot \sec}{m}$.

In experiments, the Laplace transformation method is preferred when the integrity of gas pipelines is. Laplace transformation converts equations (5) into second-order ordinary differential equations, and their general solutions will be as follows.

$$P(x,S) = \dfrac{P_i - 2aG_0 x}{S} + c_1 Sh\lambda x + c_2 Ch\lambda x - \beta G_{ut}(S)\int_{\ell_1}^{x}\delta(y-\ell_2)sh\lambda(y-x)dy \qquad (6)$$

Here, $\lambda = \sqrt{\dfrac{2as}{c^2}}$, $\beta = \sqrt{\dfrac{2ac^2}{s}}$

By applying the Laplace transform to the initial and boundary conditions and considering equation (6), we determine the coefficients $C_1$ and $C_2$. Substituting these coefficients into equation (6), we obtain the transformed equation representing the considered process, i.e., the dynamic state of the gas pipeline. These equations describe the transformed distribution of pressure along the pipeline under the given conditions.

$$P(x,S) = \frac{P_1 - 2aG_0 x}{S} + \frac{\beta G_0}{S}\frac{sh\lambda\left(x - \frac{L}{2}\right)}{ch\lambda\frac{L}{2}} + \beta G_0(S)\frac{ch\lambda(L-x)}{sh\lambda L} - \beta G_s(S)\frac{ch\lambda x}{sh\lambda L} -$$

$$- \beta G_{ut}(S)\frac{ch\lambda(L-\ell_2)ch\lambda x}{sh\lambda L} - \begin{cases} 0 \to 0 \leq x \leq \ell_2 \\ \beta G_{ut}(S)sh\lambda(\ell_2 - x) \to \ell_2 \leq x \leq L \end{cases} \quad (7)$$

We use the inverse Laplace transform rule to find the original solution of equation (7). Thus, given the initial and boundary conditions, the mathematical expression for the pressure distribution over time along the pipe length is obtained as follows:

$$P(x,t) = P_1 - 2aG_0\frac{L}{2} + 8aG_oL\sum_{n=1}^{\infty}\frac{e^{-\alpha_1 t}}{[\pi(2n-1)]^2}Cos\frac{\pi(2n-1)x}{L} +$$

$$+ \frac{2c^2}{L}\sum_{n=1}^{\infty}Cos\frac{\pi n x}{L}\int_0^t \left[G_0(\tau) - (-1)^n G_s(\tau)\right]\cdot e^{-\alpha_2(t-\tau)}d\tau + \frac{c^2}{L}\int_0^t\left[G_0(\tau) - G_s(\tau) - G_{ut}(\tau)\right]d\tau -$$

$$-\frac{2c^2}{L}\sum_{n=1}^{\infty}\int_0^t G_{ut}(\tau)\cdot e^{-\alpha_2(t-\tau)}d\tau \cdot \begin{cases} Cos\frac{\pi n x}{L}Cos\frac{\pi n(L-\ell_2)}{L} \to 0 \leq x \leq \ell_2 \\ Cos\frac{\pi n \ell_2}{L}Cos\frac{\pi n(L-x)}{L} \to \ell_2 \leq x \leq L \end{cases} \quad (8)$$

Here, $\alpha_1 = \frac{\pi^2(2n-1)^2 c^2}{2aL^2}$, $\alpha_2 = \frac{\pi^2 n^2 c^2}{2aL^2}$

## 4. Analysis and Evaluation of the Dynamics of Accidents in Gas Pipelines Using Numerical and Experimental Methodologies.

During the determination of the calculation time for transient processes during accidents in main pipeline networks, it should be considered that the simplification process does not have a significant impact, in other words, a short period is required. For example, the simplification $G_0(t) = G_s(t) = G_0$ is correct for a certain period of time, but as the duration of time increases, this equality is violated and continues until the creation of a new steady-state regime [14]. Reducing the calculation time also enhances the efficiency of dispatcher stations. This is because promptly identifying the location of a pipeline rupture allows for the activation of shut-off valves, thereby localizing the damaged section of the gas pipeline and reducing gas loss. To achieve this, it is necessary to use calculation programs based on simple principles, adopting simplification methods that have undergone theoretical and experimental tests in known areas of application. Currently, such programs are available [14].

$G_0(t) = G_s(t) = G_0$, $G_{ut}(t) = G_{ut} = constat$

When accepting these simplifications, equations (8) will be as follows:

$$P(x,t) = P_1 - 2aG_0 x + \frac{8aLG_0}{\pi^2} \sum_{n=1}^{\infty} \cos\frac{\pi n x}{L} \frac{e^{-(2n-1)^2 \alpha_3 t}}{(2n-1)^2} -$$

$$\frac{4aL}{\pi^2} G_0 \sum_{n=1}^{\infty} ((1-(-1)^n)) \cos\frac{\pi n x}{L} \frac{e^{-n^2 \alpha_3 t}}{n^2} - 2aG_{ut}\left(\frac{x^2}{2L} + \frac{\ell_2^2}{2L} + \frac{L}{3} - \ell_2\right) + \quad (9)$$

$$\frac{4aL}{\pi^2} G_{ut} \sum_{n=1}^{\infty} \cos\frac{\pi n x}{L} \cos\frac{\pi n \ell_2}{L} \frac{e^{-n^2 \alpha_3 t}}{n^2} - \begin{cases} 0 \to 0 \le x \le \ell_2 \\ 2aG_{ut}(x - \ell_2) \to \ell_2 \le x \le L \end{cases}$$

Here, $\alpha_3 = \frac{\pi^2 c^2}{2aL^2}$

Using the equations (9) mentioned above, we accept the following data to determine the lawfulness of the gas pressure variation along the length of the complex pipeline, depending on different locations of gas leakage ($\ell_2 = 0,5 \cdot 10^4$ m; $\ell_2 = 5 \cdot 10^4$ m; $\ell_2 = 9,5 \cdot 10^4$ m).

$P_1 = 55 \cdot 10^4$ Pa; $P_2 = 25 \cdot 10^4$ Pa; $G_0 = 30$ Pa·san/m; $2a = 0,1\frac{1}{\sec}$; $c = 383.3\frac{m}{\sec}$

$L = 10^5$ m: d= 0,7 m:

First, by assuming the location of the gas leakage point in the gas pipeline as $\ell_2 = 0,5 \cdot 10^4$ m, we calculate the values of the gas pressure in the pipeline every 100 seconds and every 12.5 km (including at the gas leakage point and at the 95th km), and record them in the table below (Table 1).

| x, $10^3$ m | P(x.t), $10^4$ Pa | | | | | | | | |
|---|---|---|---|---|---|---|---|---|---|
| | t=100, sec | t=200, sec | t=300, sec | t=400, sec | t=500, sec | t=600, sec | t=700, sec | t=800, sec | t=900, sec |
| 0 | 52.23 | 50.58 | 49.30 | 48.21 | 47.25 | 46.38 | 45.58 | 44.83 | 44.13 |
| 5 | 52.10 | 50.50 | 49.24 | 48.16 | 47.21 | 46.35 | 45.55 | 44.81 | 44.11 |
| 12.5 | 49.72 | 48.34 | 47.18 | 46.18 | 45.27 | 44.44 | 43.67 | 42.95 | 42.27 |
| 25 | 47.11 | 46.29 | 45.45 | 44.64 | 43.88 | 43.16 | 42.48 | 41.83 | 41.21 |
| 37.5 | 43.69 | 43.33 | 42.83 | 42.28 | 42.01 | 41.14 | 40.58 | 40.04 | 39.51 |
| 50 | 40.00 | 39.88 | 39.64 | 39.31 | 38.93 | 38.53 | 38.10 | 37.67 | 37.24 |
| 62.5 | 36.25 | 36.22 | 36.13 | 35.96 | 35.73 | 35.47 | 35.17 | 34.85 | 34.51 |
| 75 | 32.50 | 32.49 | 32.46 | 32.39 | 32.27 | 32.10 | 31.90 | 31.67 | 31.41 |
| 87.5 | 28.75 | 28.75 | 28.74 | 28.71 | 28.64 | 28.54 | 28.40 | 28.22 | 28.01 |
| 95 | 26.50 | 26.50 | 26.49 | 26.47 | 26.42 | 26.34 | 26.21 | 26.05 | 25.86 |
| 100 | 25.00 | 25.00 | 24.99 | 24.97 | 24.93 | 24.84 | 24.72 | 24.57 | 24.37 |

Table 1. The pressure variations along the length of the gas pipeline over time due to leakage. ($\ell_2 = 0,5 \cdot 10^4$ m).

Using Table 1, let's plot the graph of pressure distribution along the length of the complex gas pipeline as a function of every 200 seconds (Graph 1).

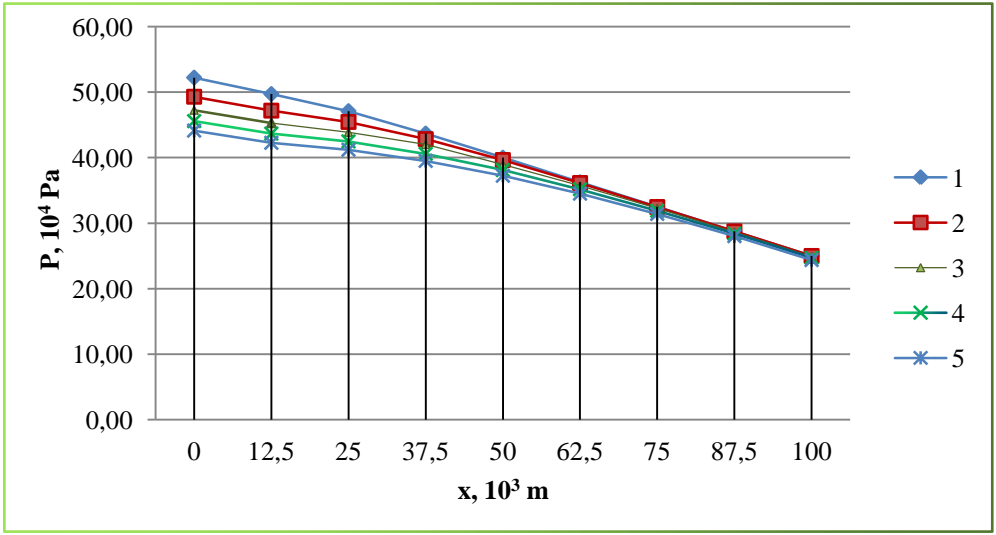

Graph 1. The distribution of pressure in the gas pipeline over time and along its length due to leakage ( $\ell_2 = 0.5 \cdot 10^4$ m, 1-t=100 sec, 2-t=300 sec, 3- t=500 sec, 4- t=700 sec, 5- t=900 sec).

Sequentially assuming the gas leakage point in the gas pipeline to be $\ell_2 = 5 \cdot 10^4$ meters, the pressure values of the gas pipeline are calculated at intervals of every 100 seconds and every 12.5 kilometers, including the 5th and 95th kilometers, and recorded in the following table (Table 2).

| x, $10^3$m | P(x,t), $10^4$Pa | | | | | | | | |
|---|---|---|---|---|---|---|---|---|---|
| | t=100, sec | t=200, sec | t=300, sec | t=400, sec | t=500, sec | t=600, sec | t=700, sec | t=800, sec | t=900, sec |
| 0 | 55 | 54.9 | 54.66 | 54.34 | 53.96 | 53.56 | 53.14 | 52.71 | 52.28 |
| 5 | 53.5 | 53.38 | 53.14 | 52.81 | 52.43 | 52.03 | 51.6 | 51.17 | 50.74 |
| 12,5 | 51.23 | 51.05 | 50.76 | 50.4 | 50 | 49.59 | 49.16 | 48.73 | 48.29 |
| 25 | 47.34 | 46.93 | 46.49 | 46.05 | 45.61 | 45.17 | 44.73 | 44.29 | 43.85 |
| 37,5 | 43.05 | 42.35 | 41.76 | 41.23 | 40.75 | 40.28 | 39.82 | 39.38 | 38.93 |
| 50 | 37.95 | 37.1 | 36.45 | 35.89 | 35.38 | 34.9 | 34.44 | 33.99 | 33.54 |
| 62,5 | 35.55 | 34.85 | 34.26 | 33.73 | 33.25 | 32.78 | 32.32 | 31.87 | 31.43 |
| 75 | 32.33 | 31.93 | 31.49 | 31.05 | 30.61 | 30.17 | 29.73 | 29.29 | 28.84 |
| 87,5 | 28.72 | 28.54 | 28.25 | 27.89 | 27.50 | 27.09 | 26.66 | 26.23 | 25.79 |
| 95 | 26.49 | 26.38 | 26.13 | 25.80 | 25.02 | 25.02 | 24.60 | 24.17 | 23.74 |
| 100 | 24.99 | 24.89 | 24.66 | 24.33 | 23.96 | 23.49 | 23.14 | 22.71 | 22.27 |

Table 2. The pressure variations along the length of the gas pipeline over time due to leakage. ( $\ell_2 = 5 \cdot 10^4$ m).

Using Table 2, let's create a graph illustrating the distribution of pressure along the length of the complex gas pipeline every 200 seconds (Graph 2).

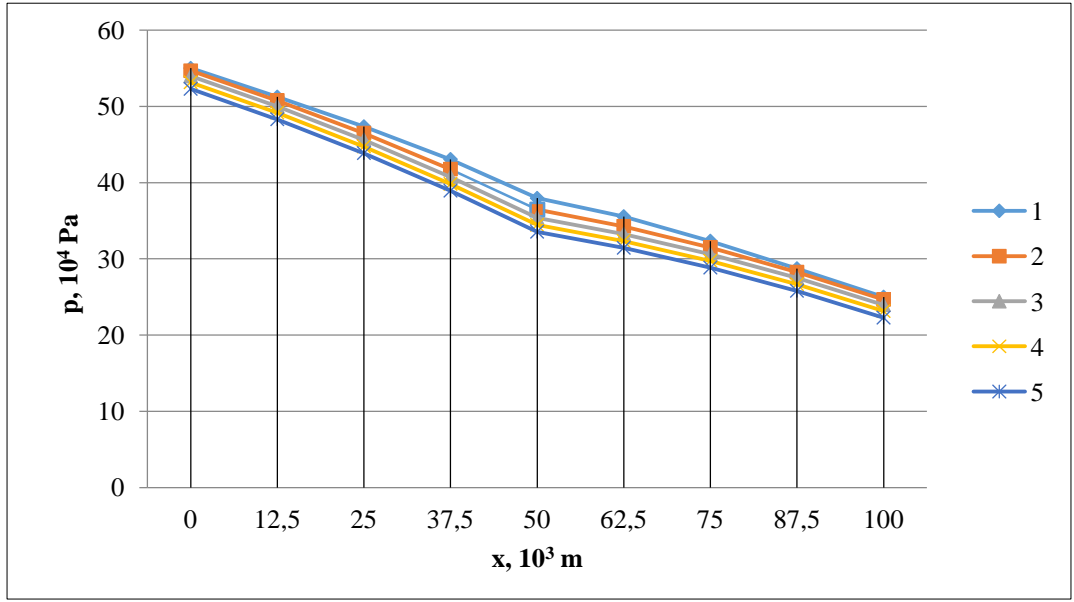

Graph 2. The distribution of pressure in the gas pipeline over time and along its length due to leakage ( $\ell_2 = 5 \cdot 10^4$ m, 1- t=100 sec, 2- t=300 sec, 3- t=500 sec, 4- t=700 sec, 5- t=900 sec)..

Finally, assuming the location of the gas leakage point at $\ell_2 = 9.5 \cdot 10^4$ m, we calculate the pressure values of the gas pipeline every 300 seconds and every 12.5 kilometers (as well as at the gas leakage point and at the 5th kilometer), and record them in the following table (Table 3).

| x, $10^3$ m | P(x,t), $10^4$ Pa | | | | | | | | |
|---|---|---|---|---|---|---|---|---|---|
| | t=100, sec | t=200, sec | t=300, sec | t=400, sec | t=500, sec | t=600, sec | t=700, sec | t=800, sec | t=900, sec |
| 0 | 55.00 | 55.00 | 54.99 | 54.97 | 54.93 | 54.84 | 54.72 | 54.57 | 54.37 |
| 5 | 53.50 | 53.50 | 53.49 | 53.47 | 53.42 | 53.34 | 53.21 | 53.05 | 52.86 |
| 12,5 | 51.25 | 51.25 | 51.24 | 51.21 | 51.14 | 51.04 | 50.9 | 50.72 | 50.51 |
| 25 | 47.50 | 47.49 | 47.46 | 47.38 | 47.26 | 47.10 | 46.90 | 46.67 | 46.41 |
| 37,5 | 43.75 | 43.72 | 43.62 | 43.45 | 43.23 | 42.96 | 42.67 | 42.35 | 42.01 |
| 50 | 39.99 | 39.88 | 39.63 | 39.30 | 38.93 | 38.52 | 38.10 | 37.67 | 37.24 |
| 62,5 | 36.18 | 35.83 | 35.32 | 34.77 | 34.20 | 33.64 | 33.08 | 32.54 | 32.00 |
| 75 | 32.10 | 31.28 | 30.44 | 29.64 | 28.88 | 28.16 | 27.48 | 26.83 | 26.21 |
| 87,5 | 27.22 | 25.83 | 24.68 | 23.67 | 22.77 | 21.94 | 21.17 | 20.45 | 19.77 |
| 95 | 23.56 | 21.96 | 20.70 | 19.62 | 18.67 | 17.81 | 17.01 | 16.27 | 15.57 |
| 100 | 22.23 | 20.58 | 19.30 | 18.21 | 17.25 | 16.38 | 15.58 | 14.83 | 14.13 |

Table 3. The pressure variations along the length of the gas pipeline over time due to leakage ( $\ell_2 = 9.5 \cdot 10^4$ m).

Using Table 3, let's create a graph illustrating the distribution of pressure along the length of the complex gas pipeline every 200 seconds (Graph 3).

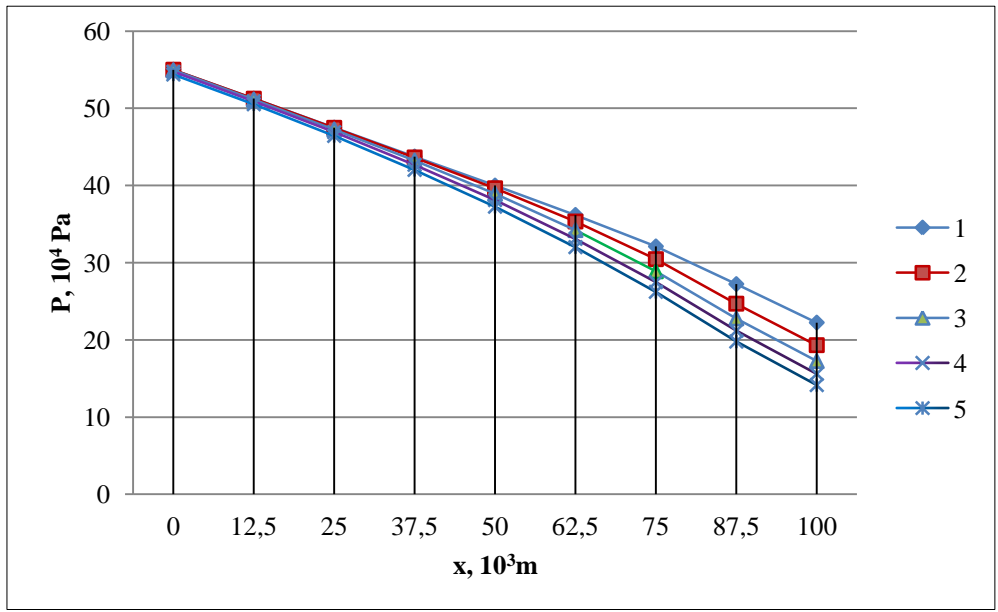

Graph 3. The distribution of pressure in the gas pipeline over time and along its length due to leakage ($\ell_2 = 9.5 \cdot 10^4$ m, 1-t=100 sec, 2-t=300 sec, 3-t=500 sec, 4-t=700 san, 5- t=900 sec).

Pressure is a key parameter in the operation of many engineering facilities, including gas supply systems. Control over this parameter is not only important during operation but also crucial during emergencies. The efficiency of gas pipelines' operation depends significantly on the control and regulation of this parameter. It should be noted that the more precise the control over this parameter, the higher the efficiency of gas pipelines' operation, which is essential for gas consumers as well. In other words, continuous control and regulation are necessary to ensure the proper operation of gas supply systems.

The analysis of the graphs and tables reveals that the quality of the longitudinal pressure variation over time varies significantly depending on the location of the leakage point. This variability is most pronounced in the pressures at the initial and final points. It appears that if an accident occurs at the beginning of the pipeline ($\ell_2 = 0.5 \cdot 10^4$ m), then at this time, the pressure drop at the beginning of the pipeline will be approximately $8.62 \times 10^4$ Pa after 600 seconds $[(P_1-P_1(0,t)] =(55-46,38) = 8,62 \cdot 10^4$ Pa. However, at the end of the pipeline, the pressure drop will be much smaller, approximately $[(P_2-P_2(L,t)] = (25-24,72) = 0,28 \cdot 10^4$ Pa, and it differs from the initial pressure drop by approximately 30 times (Graph 1 and Table 1).

If an accident (rupture) occurs at the end of the pipeline section ($\ell_2 = 9,5 \cdot 10^4$ m), then at the beginning, the pressure drop will be $0,16 \cdot 10^4$ Pa after 600 seconds $(55-54,84)=0,16 \cdot 10^4$ Pa However, at the end of the pipeline, the pressure drop will be much larger, approximately $(25-16,38)=8,62 \cdot 10^4$ Pa, which is about 53 times greater than the initial pressure drop (Graph 3 and Table 3).

Now let's consider the scenario where the accident occurs in the middle section of the pipeline ($\ell_2 = 5 \cdot 10^4$ m). In this case, the pressure drop at the beginning will be $1,44 \cdot 10^4$ Pa after 600 seconds $(55-53,56)=1,44 \cdot 10^4$ Pa. However, at the end of the pipeline, the pressure drop will also be $1,44 \cdot 10^4$ Pa after 600 seconds $(25-23,56) =1,44 \cdot 10^4$ Pa, which is equal to the initial pressure drop (Graph 2 and Table 2). The difference in pressure drops $[(P_1-P_1(0,t)]- [(P_2-P_2(L,t)]$ varies over time depending on the constant value of gas leakage. For instance, at a distance of $\ell_2 = 0,5 \cdot 10^4$ m from the leakage point, it is 2.77 at t=100 seconds, 5.67 at t=300 seconds, and

7.59 at $\tau$=500 seconds (Graph 4). This means that near the beginning of the gas pipeline, if a gas leak occurs, the value of the difference in pres increases over time, while near the end, it decreases over time. However, in the middle section of the pipeline, if a gas leak occurs, the value of the difference in pressure drops remains constant over time (Graph 4).

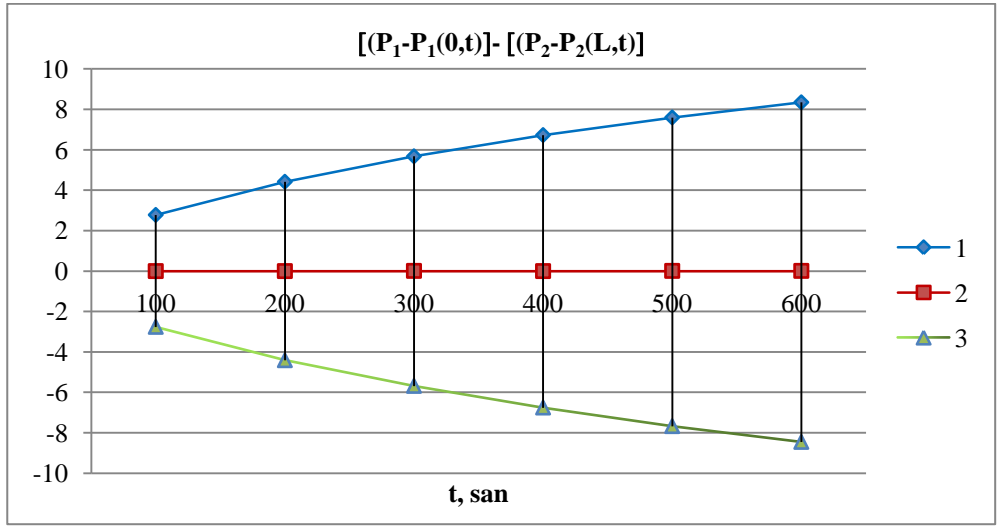

**Graph 4.** The view of the distribution of the difference in pressure drops between the beginning and end points of the gas pipeline over time and distance from the location of the pipeline's rupture (1- $\ell_2 = 0{,}5 \cdot 10^4$ m, 2- $\ell_2 = 5 \cdot 10^4$ m, 3- $\ell_2 = 9{,}5 \cdot 10^4$ m).

## 5. *Modeling of Non-Stationary Flow Parameters.*

Similarly, it can be noted that in the vicinity of the beginning of the gas leakage point of the pipeline ($\ell_2 = 0{,}5 \cdot 10^4$ m), after t=300 seconds, the value of the initial pressure drop is approximately $\left\{\dfrac{(55-49{,}30)}{(25-24{,}97)}\right\}$=190 times greater than the pressure drop generated near the end of the pipeline (Graph 1 and Table 1).

If a leakage occurs near the end of the pipeline ($\ell_2 = 9{,}5 \cdot 10^4$ m), after t = 300 seconds, the pressure drop at the beginning of the pipeline is approximately $\left\{\dfrac{(55-54{,}99)}{(25-19{,}3)}\right\}$=0,002 times smaller than the pressure drop generated at the beginning of the pipeline (Graph 3 and Table 3). If the leakage occurs in the middle section of the pipeline, then after t = 300 seconds, the pressure drop at the beginning of the pipeline equals the pressure drop generated at the end of the pipeline, $\left\{\dfrac{(55-54{,}53)}{(25-24{,}37)}\right\}$=1 (Graph 4, Table 4, and Graph 5).

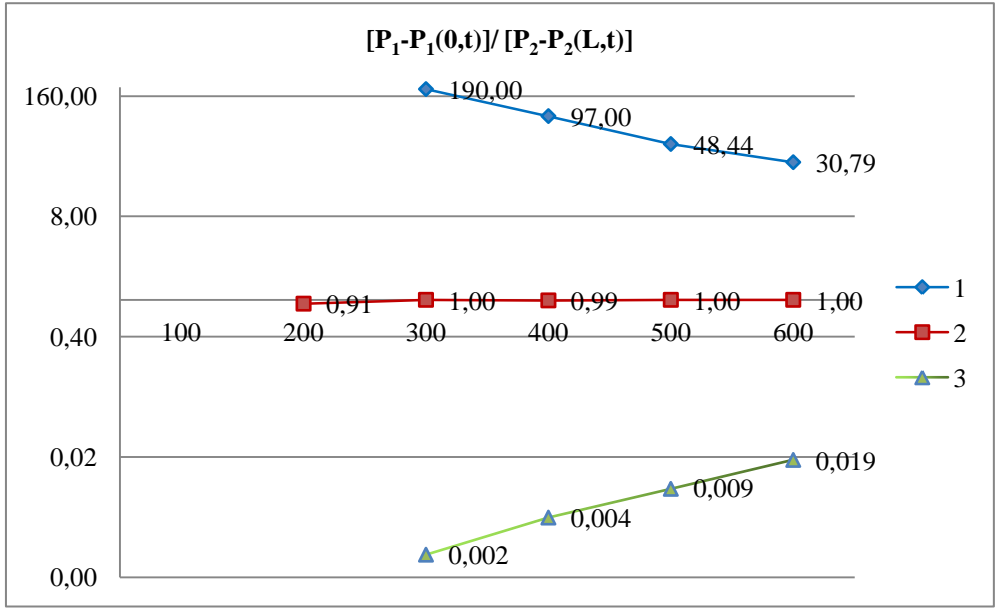

Graph 5. The view of the distribution of the ratio of pressure drops at the beginning and end points of the gas pipeline over time, depending on the location and time of the pipeline's leakage. (1- $\ell_2 = 0{,}5 \cdot 10^4$ m, 2- $\ell_2 = 5 \cdot 10^4$ m , 3- $\ell_2 = 9{,}5 \cdot 10^4$ m).

From the analysis of Graph 5, it can be concluded that the ratios of the pressure drops occurring at the beginning and end of the pipeline vary significantly depending on the location of the gas leakage point. However, as the duration of the incident increases, these ratios either decrease or increase. For instance, in the vicinity of the end of the pipeline ( $\ell_2 = 9{,}5 \cdot 10^4$ m), after t = 300 seconds, the ratio of the pressure drop at the beginning of the pipeline to that at the end is 0.002. However, at t = 400 seconds, this ratio increases to 0.004, at t = 500 seconds to 0.009, and at t = 600 seconds, it further increases to 0.019.

The further analysis of Graph 5 reveals that the value of this ratio reaches its maximum or minimum at t=300 seconds for a constant value of gas leakage, and these values vary over time. Specifically, as the duration of the incident increases, the ratio of pressure drops due to the leakage near the beginning of the pipeline decreases. However, in the middle section of the pipeline, the ratio of pressure drops due to leakage stabilizes after t=300 seconds. In subsequent times, this value remains unchanged. Thus, by optimizing the parameters of complex gas pipeline systems during incident scenarios, we determine the fixed time $t=t_1$ for non-stationary flow parameters.

## 6. Early detection method for identifying leakage points in gas pipelines.

The advantage of non-stationary gas flow models is that they allow the analytical determination of the leakage location to be correlated with the pressure change values over time (Figure 1). The difficulties in this case arise from the fact that it is impossible to completely avoid relative errors in precisely determining the exact location where fractures or cracks occur.

On the other hand, the advantage of stationary gas flow models lies in the simplicity of calculations and the accurate identification of damaged areas. However, their disadvantage is that after a pressure drop at certain points in the gas pipeline, a certain amount of time is required to establish a new stationary flow state. During this period, a significant amount of gas is released into the environment..

The occurrence of a gas leak in a gas pipeline results in the formation of a non-stationary gas flow. The duration of the non-stationary regime is determined by the transition time from one stationary state to another, and it depends on the speed of propagation of pressure fluctuations along the length of the gas pipeline. Pressure fluctuation waves in the gas flow propagate at the speed of sound. Therefore, with a gas pipeline length of 100 km, the movement duration of the wave will be approximately 5 minutes. Thus, it is not difficult to locate the leak by observing the pressure drop in the network. In this case, there is no need to wait for the formation of a new stationary regime in the gas flow.

In the fixed time interval $t=t_1$, it is possible to observe a deviation in the pressure change from its value in the operational regime. The ability to detect small gas leaks will depend entirely on the accuracy of measurement instruments. In an automated dispatch control system, calculations are made in real-time for any changes in gas pressure at selected points. In this case, the actual and calculated pressure values for each section should be compared quickly. Therefore, determining the minimum time required to deliver reliable information to the dispatch center is a crucial issue.

The result of the analysis indicates that determining the location of leakage in pipelines analytically, based on the ratio of pressure drops occurring at the beginning and end of the pipeline, is more suitable for the intended purpose. Additionally, determining the location of the leakage at the time of fixation is crucial since the ratio of pressure drops occurring at the beginning and end of the pipeline varies significantly over time from the moment of the incident. Therefore, considering the factors mentioned, such as the simplicity and accuracy of the mathematical expressions obtained as a result of optimizing non-stationary flow parameters, will ensure reliability in determining the location of leakage in complex gas pipelines through analytical methods.

Based on Graph 2, it is evident that if a leakage (incident) occurs near the middle section of the pipeline (for the specific gas pipeline considered, where $\ell_2 = 5 \cdot 10^4$ m), then the pressure drop values at the end of the pipeline over time are approximately equal to the initial values. In this case, we can assume the equality $G_0(t) = G_s(t)$.

However, based on the analysis of the pressure variation along the pipeline over time for other instances of the location of the breach. If the location of the breach occurs near the beginning of the pipeline (considering a specific gas pipeline where $\ell_2 = 0,5 \cdot 10^4$ m), then it can be observed from Graph 5 that the value of the pressure drop at the end of the pipeline remains indeterminate until $t= 300$ seconds. This is because, as you mentioned, it takes about 5 minutes for the non-stationary gas flow wave to reach and be recorded at the final point at the moment of the accident.

On the other hand, if the location of the leak occurs near the end of the pipeline (considering a specific gas pipeline where $\ell_2 = 9,5 \cdot 10^4$ m), then it also takes $t=300$ seconds for the gas flow wave to reach and be recorded at the starting point. It can be seen from Graph 5 that after $t=300$ seconds, the pressure drops at the beginning and end points of the pipeline take different values. It is understood that the change in gas consumption at the beginning and end points of the gas pipelines directly depends on the values of the pressure drops in those sections [15,16,17,18].

So, by using equations (8) and (9), we determine the following formula to identify the location of damage in non-stationary flow parameters in complex gas networks at the fixed time $t=t_1$.

$$\ell_2 = \frac{1}{2} + \varphi \cdot L \cdot \left[ \frac{1-\phi(t)}{1+\phi(t)} \right] \qquad (10)$$

Here, $\phi(t) = \dfrac{P_1 - P_1(0,t_1)}{P_2 - P_2(L,t_1)}$; $\varphi = \dfrac{2}{3} + \dfrac{2\alpha t - 4e^{-\alpha t_1} - e^{-2\alpha t_1}}{L}$ $\alpha = \dfrac{\pi^2 c^2}{2aL^2}$

As mentioned above, in accident regimes, the function $\phi(t)$ takes extreme values at $t=t1$. In this case, equation (10) is used to determine the location of the gas leakage point and the location of the operating valves. To confirm our statements and for engineering calculations, we accept the known data and parameters obtained from calculations (Tables 1-3 and Graphs 1-5).

$P_1 = 55 \cdot 10^4$ Pa; $P_2 = 25 \cdot 10^4$ Pa; $G_0 = 30$ Pa·sec/m; $2a = \dfrac{1}{san}$; $c = 383,3 \dfrac{m}{sec}$

$L = 10^5$ m: d=0.7m: $t_1 =300$ san: $\ell_2^f = 50000$ m: $\ell_2^f = 5000$ m: $\ell_2^f = 95000$ m.

Using equation (10) and considering different values of the leakage point, we determine the leakage location and compare the calculated value with the actual value. We record the results of the calculation in the following table (Table 4).

| $P_1(0,t)$, $10^4$ Pa | $P_2(L,t)$, $10^4$ Pa | t, sec | $\varphi$ | $\phi(t)$ | $\ell_2^h$, m | $\dfrac{\ell_2^f - \ell_2^h}{\ell_2^f}$ |
|---|---|---|---|---|---|---|
| | | | $\ell_2^f = 50000$ m | | | |
| 55,00 | 24,99 | 100,00 | 0,27 | 0,00 | 76915,02 | 0,54 |
| 54,90 | 24,89 | 200,00 | 0,36 | 0,91 | 51737,90 | 0,03 |
| 54,66 | 24,66 | 300,00 | 0,45 | 1,00 | 50000,00 | 0,00 |
| 54,34 | 24,33 | 400,00 | 0,53 | 0,99 | 50394,81 | 0,01 |
| 53,96 | 23,96 | 500,00 | 0,59 | 1,00 | 50000,00 | 0,00 |
| 53,56 | 23,56 | 600,00 | 0,66 | 1,00 | 50000,00 | 0,00 |
| | | | $\ell_2^f = 5000$ m | | | |
| 52,23 | 25 | 100,00 | 0,27 | - | - | - |
| 50,58 | 25 | 200,00 | 0,36 | 442,00 | 13668,81 | 0,63 |
| 49,3 | 24,97 | 300,00 | 0,45 | 190,00 | 5512,41 | 0,09 |
| 48,21 | 24,93 | 400,00 | 0,53 | 97,00 | -1437,48 | 4,48 |
| 47,25 | 24,84 | 500,00 | 0,59 | 48,44 | -6913,52 | 1,72 |
| 46,38 | 24,72 | 600,00 | 0,66 | 30,79 | -11380,80 | 1,44 |
| | | | $\ell_2^f = 50000$ m | | | |
| 55 | 22,23 | 100,00 | 0,27 | 0,000 | 76915,02 | -0,19 |
| 55 | 20,58 | 200,00 | 0,36 | 0,000 | 86495,95 | -0,09 |
| 54,99 | 19,3 | 300,00 | 0,45 | 0,002 | 94800,88 | 0,00 |
| 54,97 | 18,21 | 400,00 | 0,53 | 0,004 | 102047,13 | 0,07 |
| 54,93 | 17,25 | 500,00 | 0,59 | 0,009 | 108251,16 | 0,14 |
| 54,84 | 16,38 | 600,00 | 0,66 | 0,019 | 113114,96 | 0,19 |

Table 4. Comparison of the calculated and actual values of the gas leakage point for different leakage point values.

Table 4 shows that the calculated leakage point location is accurately determined only at $t=t_1=300$ seconds, and its relative error is minimized. The reason for the relative error lies in the simplification of the non-linear system of equations and the simplification of mathematical expressions.

However, it is accepted for engineering calculations and dispatcher stations. Indeed, the principles obtained through the calculation of mathematical expressions derived from the study of physical processes and the application of the new calculation scheme are suitable for practical use in experience.

## 7. Conclusion

This study presents a novel approach for detecting gas leaks in transmission pipelines by analyzing unsteady flow parameters and optimizing pressure monitoring techniques. The key findings and contributions are summarized as follows:

A refined mathematical model for unsteady gas flow in pipelines has been developed to enhance the accuracy of pressure distribution analysis in leakage scenarios. The study introduces a computational framework that improves existing gas leakage detection methods.

A new calculation scheme has been formulated to determine the minimum fixation time ($t = t_1$) required for reliable leak detection. The analytical method proposed enables real-time detection of gas leaks, enhancing the reliability of automated pipeline control systems. The findings contribute to the design of safer and more efficient gas transportation networks.

To further enhance the effectiveness and applicability of the proposed method, future research should:

Develop hybrid detection models integrating pressure, temperature, and acoustic signals for a more comprehensive analysis.

Investigate the impact of external environmental conditions, such as soil properties and temperature fluctuations, on gas leak propagation.

Apply machine learning and artificial intelligence techniques to improve real-time detection accuracy and reduce false alarms.

Conduct large-scale experimental validation under diverse operational conditions to refine the model and optimize its practical implementation.

By addressing these aspects, future advancements can enhance the robustness and scalability of gas leakage detection systems, ensuring greater safety and efficiency in gas transportation networks.

## References

[1]. A.T.Bakesheva, (2019) "Analysis of reliable operation of main gas pipelines". In Proceedings of the IV International Scientific-Practical Conference, Membership in the WTO: Prospects of Scientific Researches and International Technology Market, Volume 1, Vancouver, Canada.
[2]. F.Matiko, L. Lesovoy, V. Dzhigirei, (2016) "Improvement of mathematical models of natural gas flow during its outflow from a damaged gas pipeline". Bulletin of the Engineering Academy of Ukraine (1), pp.224–230 .
[3]. K. Brünenberg, D. Vogt, M. Ihring, (2004) "Additional functionalities of model-based leak detection systems to improve pipeline safety and efficiency," Pipeline Technology Journal, no. 1, pp. 38-44, .
[4]. L.V. Lesovoy, L.V. Blyzniak, (2020)"Determination of natural gas pressure at the points of its removals in a gas pipeline with branches". Quality control methods and devices № 12, pp.88-91 .
[5]. T.I. Lapteva, M.N. Mansurov, (2006) "Detection of leaks during unsteady flow in pipes". Oil and Gas Business, No. 2, : http://ogbus.ru/article/obnaruzhenie-utechek-prineustanovivshemsya-techenii-v-trubax/.
[6]. U. J.Obibuike, A. Kerunwa, M. Udechukwu, R. C. Eluagu, A. C. Igbojionu, S. T. Ekwueme, (2020) "Mathematical Approach to Determination of the Pressure at the Point of Leak in Natural Gas Pipeline". Int. Journal of Oil, Gas and Coal Engineering Vol.8(1), pp.22–27.
[7]. V. Dzhyhyrei , F. Matiko, (2024)" Investigation of Changes in Natural Gas Parameters along a Damaged Gas Pipeline". Energy Engineering And Control Systems, Volume 10, Number 1, pp. 64- 71.

[8]. Y. V. Doroshenko, (2020) "Modeling of gas leaks from gas pipelines in emergency situations", Visnyk VPI, Vol.3, c. 22–28, Worm.

[9]. I.G. Aliyev, (2024) "Methodology for selecting conditions of non-stationary gas pipeline exploitation during the reconstruction stage", DOI: 10.31618/ESU.2413-9335.2024.1.119, Eastern European Scientific Journal, Том 1, № 04, Pages 3-37.

[10]. S.I.Rakhmatullin, A.G.Gumerov, V.V. Vanfatova, (2001) "Assessment of the dynamic volume balance of oil in a pipeline with gravity sections". Pipeline Transport of Oil, Supplementary Issue, No. 3, pp.24–27.

[11]. T.L.Council, D.J. Honey, M.L Cox, (2000) "Environmental solutions - key successful south Texas line installation". Oil & Gas Journal, 98(30), pp.70–72 .

[12]. M.G.Sukharev, R.V. Samoylov, (2016)" Analysis and control of stationary and non-stationary gas transport regimes". Moscow: Gubkin Russian State University of Oil and Gas .

[13]. M.V. Lurye, (2017) "Theoretical foundations of pipeline transport of oil, petroleum products, and gas". Moscow. [14]. V. Zapukhliak, L. Poberezhny, P. Maruschak, V. Grudz , R. Stasiuk, J.Brezinová, A. Guzanová, (2019) "Mathematical Modeling of Unsteady Gas Transmission System Operating Conditions under Insufficient Loading". Energies Vol.12 (7), pp.1325, https://doi.org/10.3390/ en12071325

[15]. A.G.Vanchin, (2014) "Methods for calculating the operation mode of complex main gas pipelines". Oil and Gas Business, No. 4, pp.192–214.: http://ogbus.ru/issues/4_2014/ogbus_ 4_2014_p192-214_ VanchinAG_ru.pdf.

[16]. E. Guelpa, A. Bischi, V. Verda, (2019) "Towards future infrastructures for sustainable multi-energy systems": a review. Energy, 184, pp. 2–21.

[17]. I.G. Aliyev, (2024) "Analysis Of The Complex Gas Pipeline Exploitation Process In Various Operating Modes". Journal of Theoretical and Applied Mechanics, Vol. 54, Issue 4, pp. 446 - 464.

[18]. P.A. Velmisov, A.V.Gladun, (2016) "On the management of pipeline dynamics". Journal of the Middle Volga Mathematical Society, Vol.18(4), pp. 89–97.